\documentclass[12pt]{amsart}
\usepackage{amssymb}
\usepackage{graphicx}

\newcommand{\dt}[1]{\mbox{det}({#1})}

\newcommand{\R}{\mathbb R}

\newtheorem{thm}{Theorem}

\newtheorem{lemma}{Lemma}

\newcommand{\Pf}{\noindent {\bf Proof. }}

\begin{document}
\title[Invisible knots and rainbow rings]{Invisible knots and rainbow rings: knots not determined by their determinants}

\author[J.~Godzik]{James Godzik}
\address{
Department of Mathematics,
UC, Berkeley,
Berkeley, CA 94720-3840}
\author[N.~Ho]{Nancy Ho}
\address{Tapestry Solutions}
\author[J.~Jones]{Jennifer Jones}
\address{
Department of Mathematics,
Colorado State University,
Fort Collins, CO 80523-1874}
\author[T.W.~Mattman]{Thomas W.~Mattman}
\address{
Department of Mathematics and Statistics,
California State University, Chico,
Chico, CA 95929-0525}
\author[D.~Sours]{Dan Sours}
\address{
Chico High School,
Chico, CA 95926}

\begin{abstract}
We determine p-colorability of the paradromic rings. These rings arise by generalizing the well-known experiment of bisecting a Mobius strip. Instead of joining the ends with a single half twist, use $m$ twists, and, rather than bisecting ($n = 2$), cut the strip into $n$ sections. We call the resulting collection of thin strips $P(m,n)$. By replacing each thin strip with its midline, we think of $P(m,n)$ as a link, that is, a collection of circles in space. Using the notion of $p$-colorability from knot theory, we determine, for each $m$ and $n$, which primes $p$ can be used to color $P(m,n)$.

Amazingly, almost all admit 0, 1, or an infinite number of prime colorings! This is reminiscent of solutions sets in linear algebra. Indeed, the problem quickly turns into a study of the eigenvalues of a large, nearly diagonal matrix.

Our paper combines this explicit calculation in linear algebra with a survey of several ideas from knot theory including colorability and torus links.
\end{abstract}

\maketitle

\begin{figure}[p]
\begin{center}
\includegraphics[scale=.40]{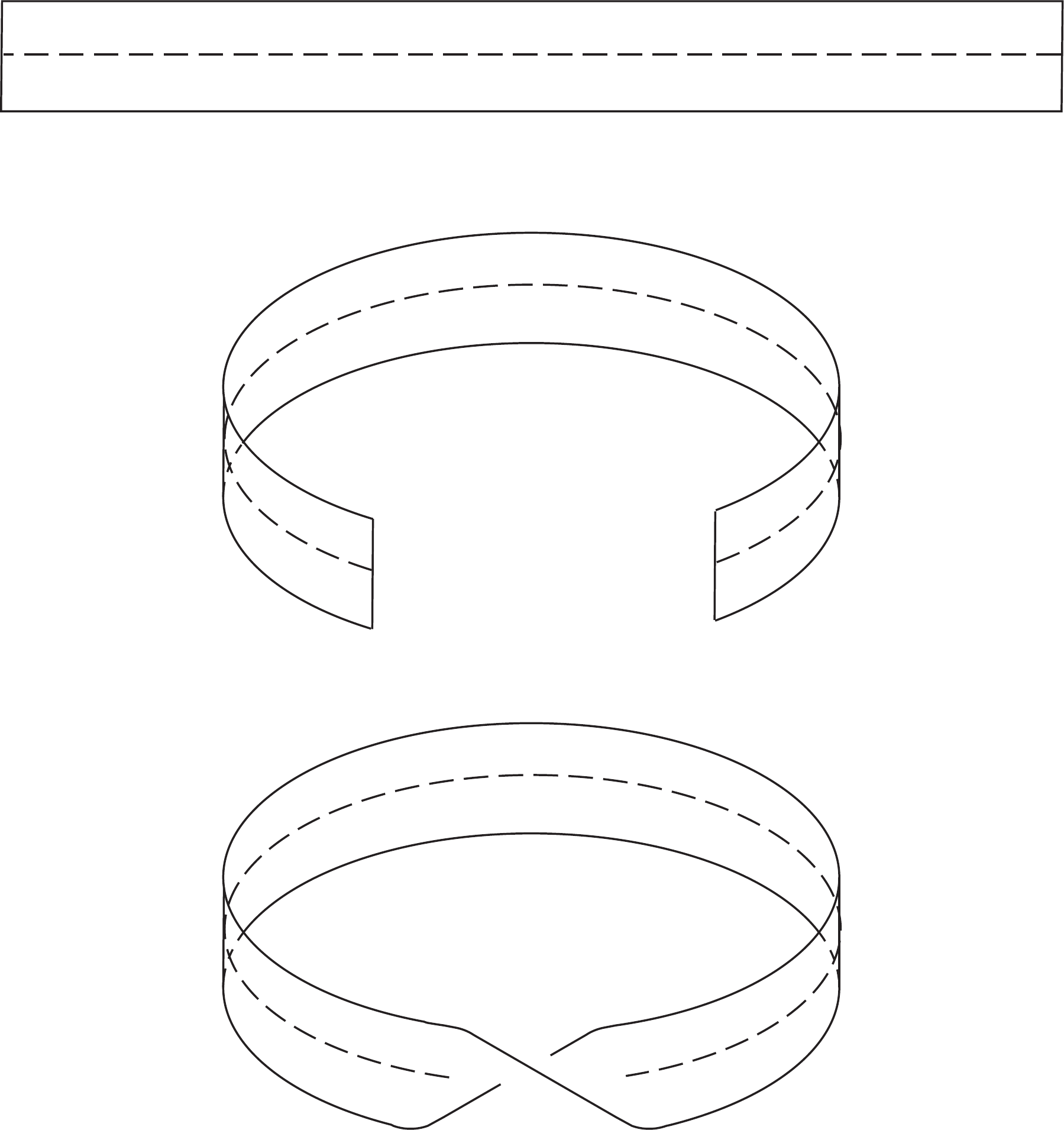}
\caption{\label{figbismob}%
Bisecting a M\"obius strip. After joining the ends with a half-twist, cut along the dashed line. What is the result?}
\end{center}
\end{figure}

M\"obius strip experiments are surefire triggers of Aha! 
experiences, even in very young audiences. Maybe
you don't remember the first time someone challenged you to color one
side blue and the other red, or asked you to guess the result of 
cutting a M\"obius strip in half, but you surely recall the outcome.
(If not, we encourage 
you to put aside the magazine for a moment, gather up some paper, tape, and scissors, and remind yourself what a bisected M\"obius strip looks like. See Figure~\ref{figbismob}).

\textbf{FIGURE 1 GOES NEAR HERE}.

As part of a research experience for undergraduates (REU), 
we discovered that generalizing these experiments results in
many more confounding constructions. Rather than simply bisecting
the M\"obius strip, try cutting it into $n$ sections. Or, instead of joining the ends of the strip with a single half twist, make two twists, or three, or, in general, $m$ half twists.  You have just created examples of {\em paradromic rings}, which
we'll denote $P(m,n)$. (We first learned of these constructions from the delightful book of Ball and Coxeter~\cite{BC}.)

\textbf{FIGURE 2 GOES NEAR HERE}.

\begin{figure}[p]
\begin{center}
\includegraphics[scale=.50]{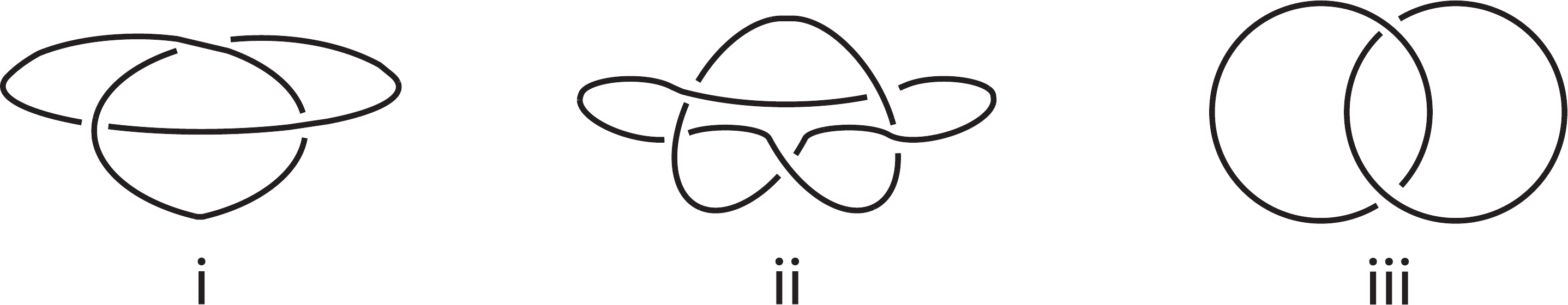}
\caption{\label{figPn2}%
Some paradromic rings with $n = 2$ (bisection)
i) $P(3,2)$, the trefoil knot, ii) $P(5,2)$, the pentafoil knot, and iii)
$P(2,2)$, the Hopf link. }
\end{center}
\end{figure}

Figure~\ref{figPn2} shows some of the results.  
Now that you have your scissors out (Get them!),
you'll find that $P(2,2)$ (bisect a strip after making a full twist) 
gives two strips of paper linked as in a chain. When $m$ is odd 
(an odd number of half twists), bisection results in a single strip, albeit knotted up. 

Having generated a nice pile of shredded strips, you'll start to wonder,
``How can we organize this tangled mess?" The very language
we are using suggests knot theory as the appropriate setting.
A {\em knot} is a simple closed curve in space, like $P(3,2)$ or $P(5,2)$ of 
Figure~\ref{figPn2}, whereas a {\em link}, like $P(2,2)$, is a collection of 
such embedded circles, called the {\em components} of the link. 
A knot, then, is a link of one compontent, and we'll use the phrase 
`links that are not knots' for those having two or more closed curves.
To realize the paradromic rings as curves,
replace each strip with its midline (or, equivalently, shrink the width of the strip to zero). 

Somehow forgetting all about the challenges of coloring 
M\"obius strips, the REU team set out to color these curves.
This is akin to edge-coloring of graphs. Just as each graph has a chromatic number,  the determinant of link $L$, $\dt{L}$, characterizes its
colorability. We'll explain how
to calculate this non-negative integer later. 
For now, it's enough to know that
$L$ is {\em $p$-colorable} if the prime $p$ divides $\dt{L}$. 
In this paper we organize the paradromic rings by colorability.
For each $m$ and $n$, we will determine the primes $p$ for which 
$P(m,n)$ is $p$-colorable.

If the word `determinant' makes you smile, you're in luck. 
In the REU, we were surprised by how quickly 
this problem in knot theory turned into a
cute exercise in linear algebra. Rather than 
calculating determinants, we'll investigate the eigenvalues of a
large, nearly diagonal matrix. There'll be some proof by pictures too, but the 
essence of our argument is algebraic.

The real Aha!, however, came when we understood that, much like the 
M\"obius strip, the paradromic rings resist coloring. Most of the 
knots in this family have determinant equal to one. 
This means they are not colorable for any prime ({\bf no solutions}). 
We call them {\em invisible knots}, following Butler et al.~\cite{BCD}.
Links of more than  one component have even determinant, and are, therefore, not invisible. Still, these paradromic rings that are not knots valiantly defy us as best they can given this constraint. 
Many have determinants that are a power of two. These we call {\em nearly invisible} as they can be colored only by the prime $p = 2$ ({\bf one solution}). 
So long as $n \neq 2,4$, the remaining paradromic rings have
$\dt{P(m,n)} = 0$. We refer to such links as {\em rainbow rings} as they can be colored by every 
prime ({\bf infinite solution set}).

In the end, the determinant is not very discriminating in separating out the 
paradromic rings. With a few exceptions, it partitions this doubly infinite family 
into only three different classes. Moreover, these classes turn out to be pathological, admitting either zero, one, or an infinite number of prime
colorings. On the other hand, perhaps this type of outcome is 
exactly what you would expect from what is, ultimately, a problem in linear algebra.

We've organized our paper as follows. 
In the next section we explain the notion of $p$-colorability of a link. 
In Section 2 we show that the paradromic rings fall into two families. If $mn$ is even, then 
we can arrange
$P(m,n)$ on the surface of a torus; it is a {\em torus link}. If $mn$ is odd, then $P(m,n)$ is a torus link with the addition of 
a circle that follows the core of the torus. In the third section we 
use linear algebra to analyze the colorability of the paradromic rings. 
The knots $P(m,1)$ are invisible, so we can assume $n > 1$.
When $mn$ is even, $P(m,n)$ is a rainbow ring except for two cases: 1) when $n = 2$ or $4$; and 2) when $n$ and $m/2$ are both odd (in which case it's nearly invisible). When $mn$ is odd, $P(m,n)$ is nearly invisible.

\section{Coloring Links}

While the determinant is convenient for organizing our 
results and defining invisible knots and rainbow rings, 
we will not calculate $\dt{P(m,n)}$ explicitly.
Rather, we define $p$-colorabiliy using link diagrams. 
A {\em diagram} is a projection
of the link into the plane with gaps left in the curve to show where it crosses over itself. For example, Figure~\ref{figPn2} consists of diagrams of the links $P(3,2)$, $P(5,2)$, and $P(2,2)$.

\textbf{FIGURE 3 GOES NEAR HERE}.

\begin{figure}[p]
\begin{center}
\includegraphics[scale=.80]{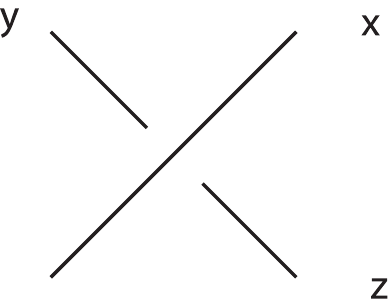}
\caption{\label{figpcross}%
Arcs are colored so that, at crossings, $2x \equiv y + z \pmod{p}$. 
The arc labeled $x$ is called an overarc, and $y$ and $z$ are underarcs.}
\end{center}
\end{figure}

\begin{figure}[p]
\begin{center}
\includegraphics[scale=.60]{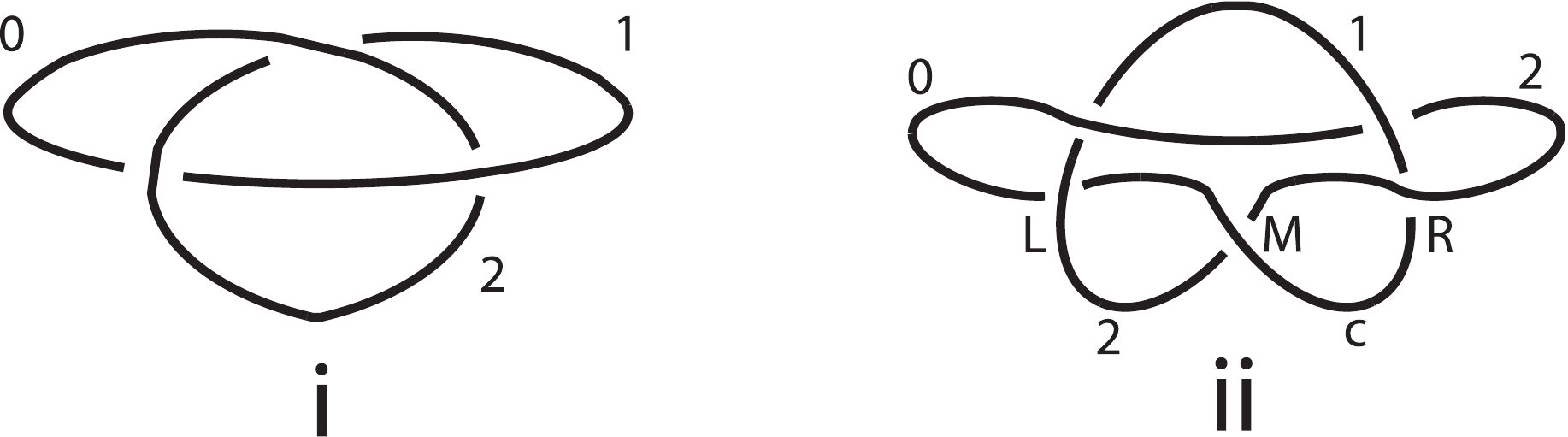}
\caption{\label{figtri}%
i) The trefoil can be tricolored.
ii) There is no way to choose a color $c$.
}
\end{center}
\end{figure}

Given a prime $p$, a diagram of a link is {\em $p$-colorable} if we can label its arcs with {\em colors} chosen from $0$, $1, \ldots,  p-1$ such that 
\begin{enumerate}
\item more than one color is used, and 
\item at each crossing the
colors satisfy the equation $$2x \equiv y+z \pmod{p}$$
\end{enumerate}
(see Figure~\ref
{figpcross}). 
A link is $p$-colorable if it has a $p$-colorable diagram. For example, Figure~\ref{figtri}i shows 
that the trefoil knot is $3$-colorable.

Condition 1 rules out the trivial solution where
every arc has the same color. Whatever the link and whatever
the prime $p$, if all arcs have color 1 (for example), condition 2 will 
hold at every crossing. Without condition 1, every link would be colorable for every $p$. You can think of the second condition 
as balancing the colors on the overarc with those on the underarcs.  There
are four lines radiating from the center of the crossing, the two on top each carrying an $x$ and the ones on the bottom carrying a $y$ and 
a $z$. Condition 2 equates the two $x$'s on top with 
the $y$ and $z$ below.

Condition 2 has a particularly nice interpretation in the case of tricolorability,
when $p = 3$. A little thought will convince you that $2x \equiv y+z \pmod{3}$ implies
either $x = y = z$ or else $\{x,y,z\} = \{0,1,2\}$. 
A link is {\em tricolorable}, then, if you can label its arcs with $0, 1, 2$ such that at least two colors are used and, at each crossing, either exactly one color, or else all three colors, appear.

\textbf{FIGURE 4 GOES NEAR HERE}.

We've mentioned that the trefoil knot $P(3,2)$ is tricolorable (Figure~\ref{figtri}i);  
let's see why the pentafoil $P(5,2)$ is not.
In Figure~\ref{figtri}ii, in trying to tricolor this knot, we have labeled four of its five arcs. 
All three colors appear at both of the top crossings,
which is consistent with condition 2. 
It's impossible, however, to assign a color $c$ to the remaining arc.
That arc is part of three crossings, one at left (L), one at right (R), and one in the middle (M). At the left crossing, the other arcs already carry $0$ and $2$, so condition 2 forces $c = 1$. On the other hand, the crossing at right obliges $c = 0$ since $1$ and $2$ already appear there. This shows that 
there is no consistent way to choose the color $c$.  Note that the middle crossing implies $c = 2$ because there are already two color $2$ arcs at that crossing.

To complete the argument that the pentafoil is not tricolorable, 
see if you can show that, no matter how the first four arcs are colored, it is impossible to choose a color $c$ for the final arc.  
(Hint: By symmetry, you may assume
the left arc is colored $0$ as in Figure~\ref{figtri}ii. There are three choices 
for the color of the top arc. With those two arcs labeled, condition 2 determines the color of two other
arcs. In other words, up to symmetry, there are only three legitimate ways to color the first four arcs.)

When $p=2$, condition 2 becomes
$y \equiv z$. At each crossing, the two 
underarcs must have the same color.  
Each component of the link, then, will be all of one 
color. As condition 1 requires we use both colors, 
a link will be $2$-colorable exactly if it has at least 
two components. 
As mentioned in the introduction, 
we say a link is nearly invisible if $p=2$ is the only coloring.

We want to use $p$-colorability to organize the paradromic rings.
It's an invariant of links, which means if a diagram admits a $p$-coloring for 
a given $p$, then any equivalent link will also have a $p$-colorable diagram.
In knot theory, we consider two links equivalent if there's
a way to move one around in space to look just like the other without 
ever having to pass the curve through itself. 
For a more precise description of link equivalence and the cute 
proof that $p$-coloring is an invariant, we recommend
Adams's {\em The Knot Book}~\cite{A} or Livingston's {\em Knot Theory}~\cite{L}.

\textbf{FIGURE 5 GOES NEAR HERE}.

\begin{figure}[p]
\begin{center}
\includegraphics[scale=.45]{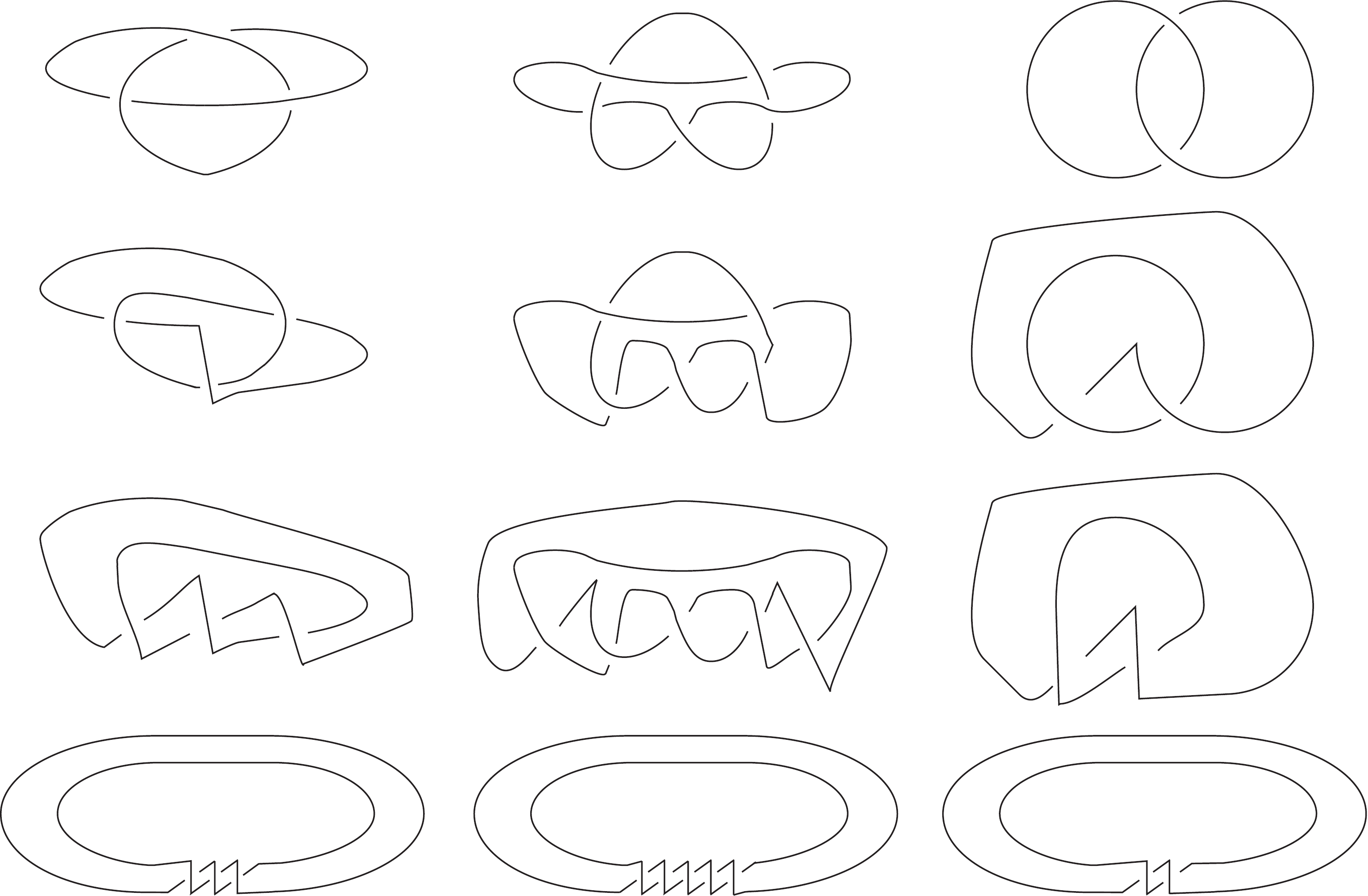}
\caption{\label{FigP352braid}%
We can redraw $P(3,2)$, $P(5,2)$, and $P(2,2)$ as at the bottom of the figure.}
\end{center}
\end{figure}

Each column of
Figure~\ref{FigP352braid} consists of four diagrams
of the same link. We've shown how the knot at left, 
$P(3,2)$, is  $3$-colorable using the top diagram. This means the three diagrams below it are also $3$-colorable, as you can easily confirm.
On the other hand, we've argued that the knot represented in the middle column, 
$P(5,2)$, is not $3$-colorable. Since $p$-colorability is a link invariant, $P(3,2)$ and $P(5,2)$ are not equivalent. There's no way to move any knot in the $P(5,2)$ column around in space to make it look just like one in the $P(3,2)$ column.
See if you can show that the third link in 
the figure, $P(2,2)$, is different from the first two. (Hint: try $5$- and $2$-colorings.
How are the $p$-colorings of  $P(m,2)$ determined by $m$?)

If you've been impatient for the linear algebra, your wait is over. But first a spoiler alert. If you haven't had a chance to see how $P(2,2)$ differs in colorability 
from the other two links in Figure~\ref{FigP352braid}, you really ought to try it before 
reading on. Remember $2$-coloring is easy. A link is $2$-colorable exactly if 
it has at least two components. You should also investigate which links in Figure~\ref{FigP352braid} are $5$-colorable.

\textbf{FIGURE 6 GOES NEAR HERE}.

\begin{figure}[p]
\begin{center}
\includegraphics[scale=.60]{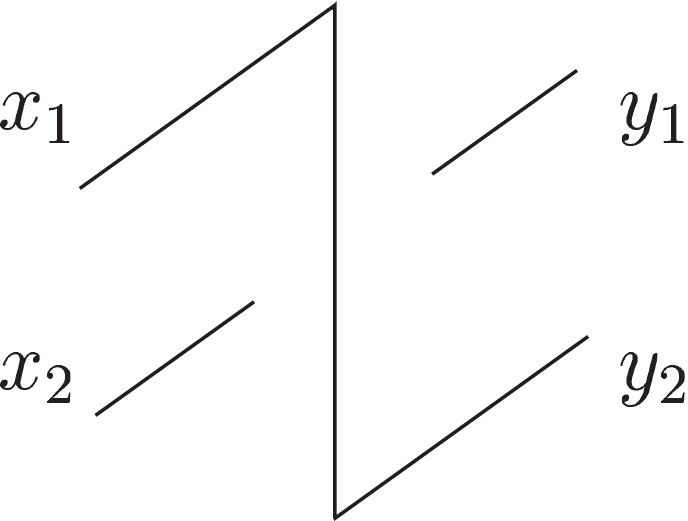}
\caption{\label{FigPm2pattern}%
Repeat this pattern $m$ times to form a $P(m,2)$ link.}
\end{center}
\end{figure}

We will now use linear algebra to prove that
$P(m,2)$ is $p$-colorable if and only if $p$ divides $m$.
The key observation is suggested by  Figure~\ref{FigP352braid}. 
To build link $P(m,2)$, repeat the Figure~\ref{FigPm2pattern} pattern
$m$ times and then join up the loose ends.
Use $x = (x_1, x_2)$ 
to color the arcs entering Figure~\ref{FigPm2pattern} at left. 
Then the arcs leaving at right are $y = (y_1, y_2)$ where $y_2 = x_1$ and condition 2 tells us that $y_1 \equiv 2x_1 - x_2 \pmod{p}$. In other words, $y \equiv T x \pmod{p}$ where 
$T= \left(
\begin{array}{rr}
2 & -1 \\
1 & 0
\end{array}
\right)
$.

\textbf{FIGURE 7 GOES NEAR HERE}.

\begin{figure}[p]
\begin{center}
\includegraphics[scale=.60]{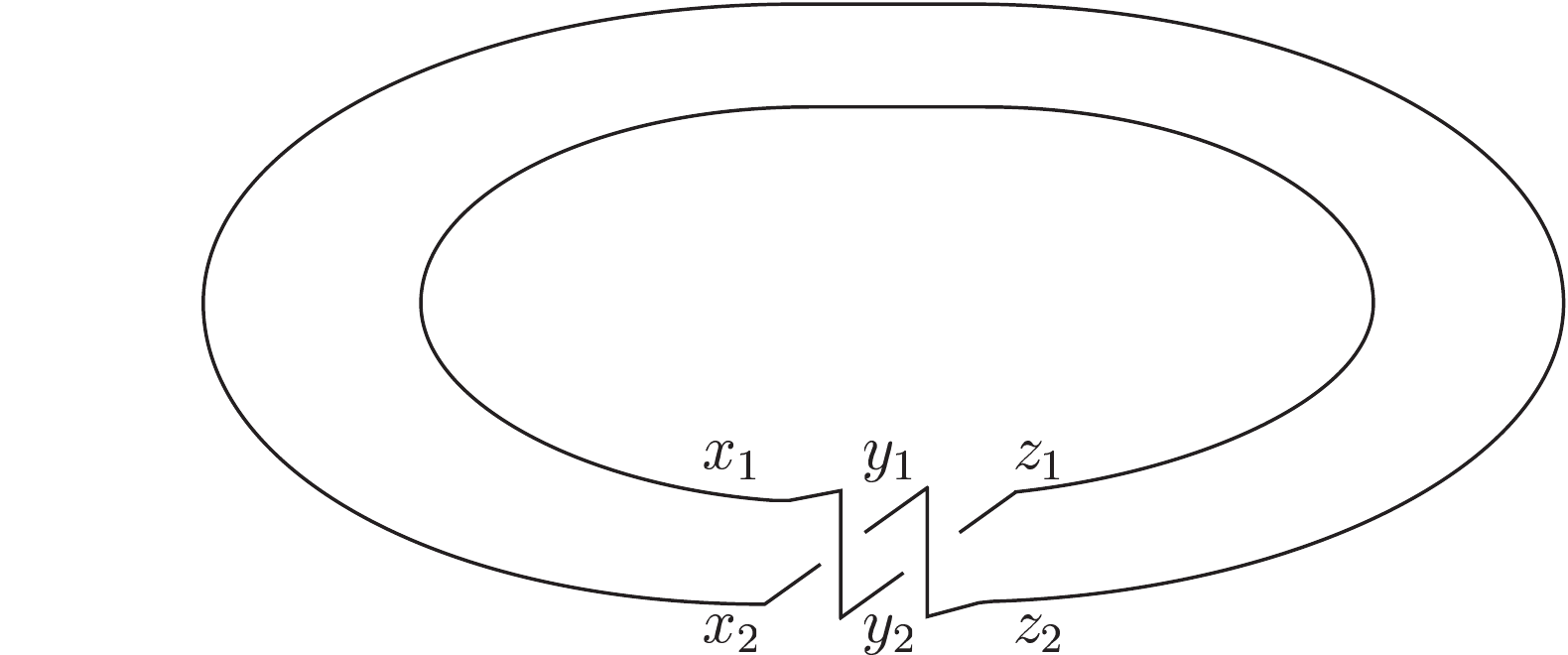}
\caption{\label{FigHopf}%
The Hopf link $P(2,2)$.}
\end{center}
\end{figure}

For the Hopf link, $P(2,2)$ (Figure~\ref{FigHopf}), we repeat the pattern two times. 
Beginning with arcs labeled $x$ at left, after going through the
pattern once, we'll have colors $y$ where $y  \equiv T x$. Passing through the pattern a second time, we have colors $z \equiv T y \equiv 
T^2 x$. Notice that by going around the top of the link these $z$ arcs at right are identified with the $x$ arcs we started with on the left.  
In other words, $x = z \equiv T^2 x$. Thus, 
$x$ represents a coloring of the Hopf link if
$x \equiv T^2 x$. 

In general, for $P(m,2)$, we pass through the Figure~\ref{FigPm2pattern} pattern $m$ times. See Figure~\ref{FigP352braid} for examples with $m = 3,5,2$.
This means a valid coloring requires $x \equiv T^m x$. Equivalently, $x$ must
satisfy the eigenvector equation: $(T^m - I) x \equiv 0$.

For any color $c$, we call $x = (c,c)$ a {\em constant vector}.
Then, $Tx = x$, so constant vectors solve the eigenvector equation.
But this means we've colored every arc $c$, violating
condition 1. Thus, $p$-colorings of $P(m,2)$ 
correspond to non--constant $\lambda = 1$ eigenvectors of $T^m$
mod $p$.

Using induction, we find $T^m - I  = \left(
\begin{array}{rr}
m & -m \\
m & -m
\end{array}
\right).$
As we mentioned, vectors of the form $(c,c)$ 
are in the null space of this matrix.  The link $P(m,2)$ will be $p$-colorable exactly when there is some other, non-constant vector in the mod $p$ null space of $T^m - I$. That means
the null space is two-dimensional so that the matrix is in fact the 
zero matrix mod $p$. Therefore, the link $P(m,2)$ is $p$-colorable if and only if $p$ divides $m$.

In Section~3, we will use this approach to determine the $p$-colorability
of the paradromic rings.

\section{Paradromic rings and torus links}

\textbf{FIGURE 8 GOES NEAR HERE}.

\begin{figure}[p]
\begin{center}
\includegraphics[scale=.5]{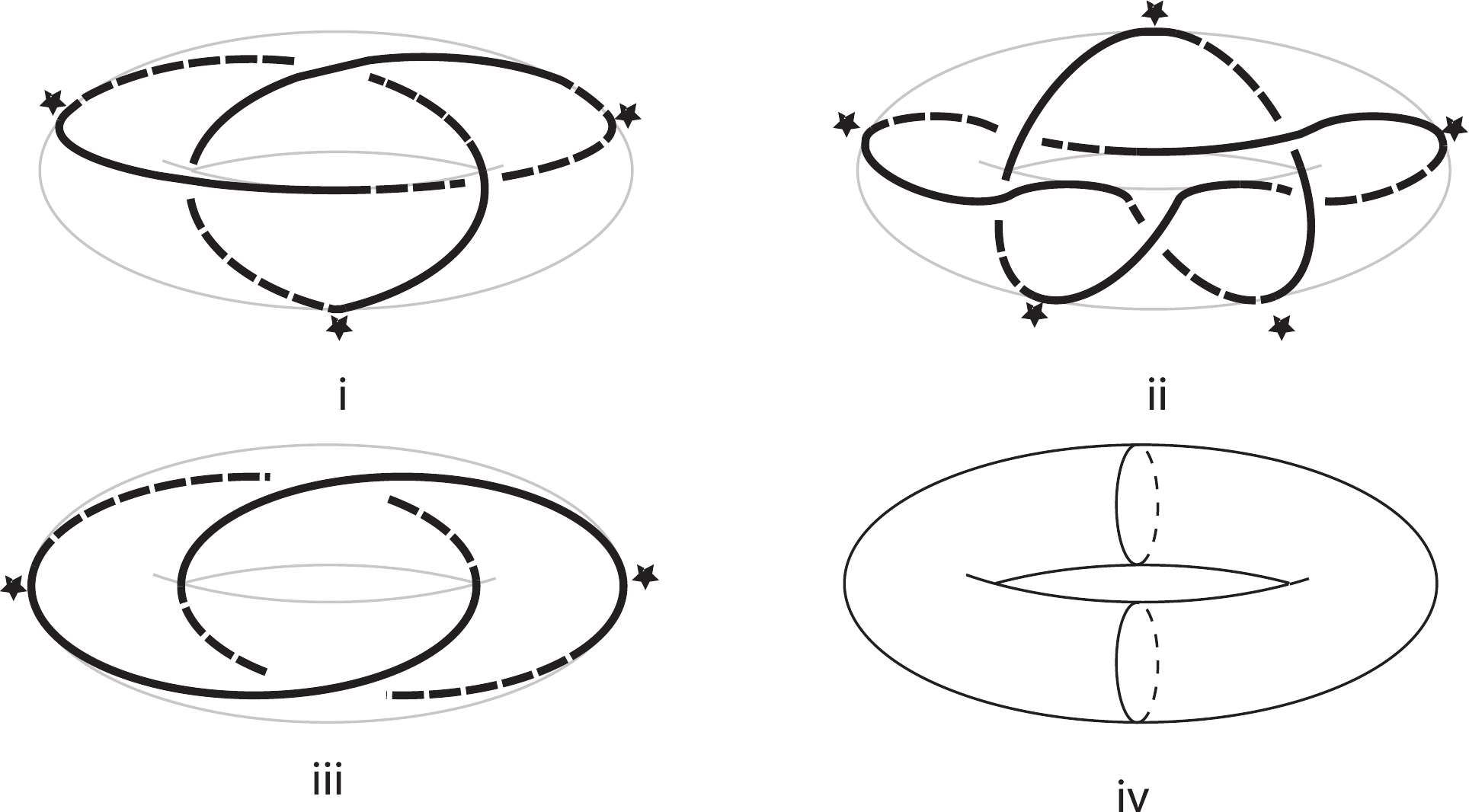}
\caption{\label{figtorex}%
The i) trefoil knot, ii) pentafoil knot, and iii) Hopf link are torus links as they can be made to lie on a torus (the surface of a doughnut, see iv). Dashed lines represent parts of the curve on the far side of the torus. }
\end{center}
\end{figure}

Paradromic rings enjoy a close connection with torus links that we will exploit
to understand their $p$-colorability.
Figure~\ref{figtorex} shows how the trefoil knot, pentafoil knot, and Hopf link
are torus links, meaning we can realize them as curves that lie flat on a torus. 
This is similar to defining a planar graph as one we can put
in the plane with no edges crossing. Links that lie in the plane are
called trivial links; they're simply collections of disjoint circles with no crossings whatsoever. 
The torus links, in contrast, are  an important family that have long
intrigued knot theorists.

\textbf{FIGURE 9 GOES NEAR HERE}.

\begin{figure}[p]
\begin{center}
\includegraphics[scale=.65]{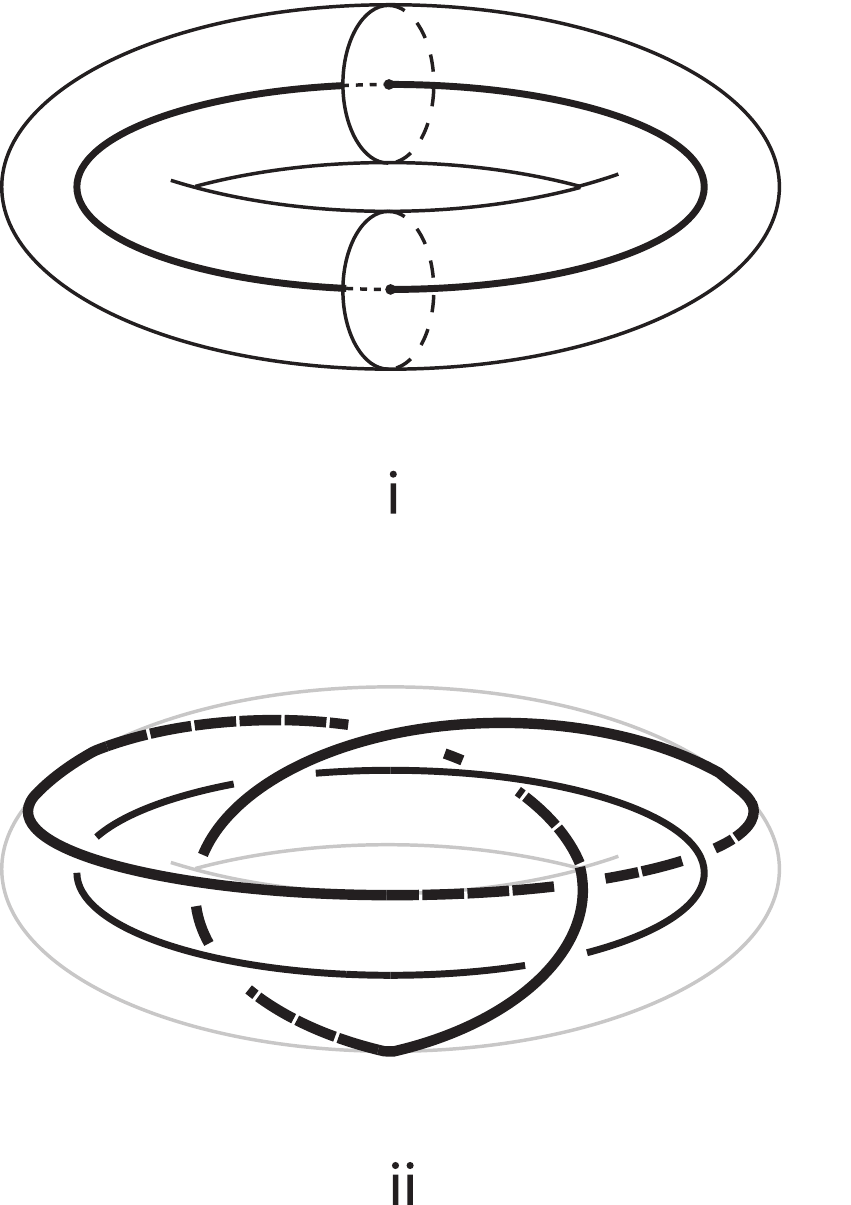}
\caption{\label{figcore}%
i) The core of the torus meets every cross-sectional disk
in its center. ii) $P(3,3)$ consists of a trefoil knot that lies on the
torus 
along with a second component along
the core of the torus.}
\end{center}
\end{figure}

We will show that each $P(m,n)$ is either a  torus link or else a torus link 
together with an additional component
that follows the {\em core} of the torus (see Figure~\ref{figcore}i).
The core is a curve inside the torus that intersects every cross-sectional disk at its center. For example, Figure~\ref{figcore}ii shows
that $P(3,3)$ consists of two components: 
the trefoil, which is a torus knot (compare Figure~\ref{figtorex}i), and the core.

\textbf{FIGURE 10 GOES NEAR HERE}.

\begin{figure}[p]
\begin{center}
\includegraphics[scale=.50]{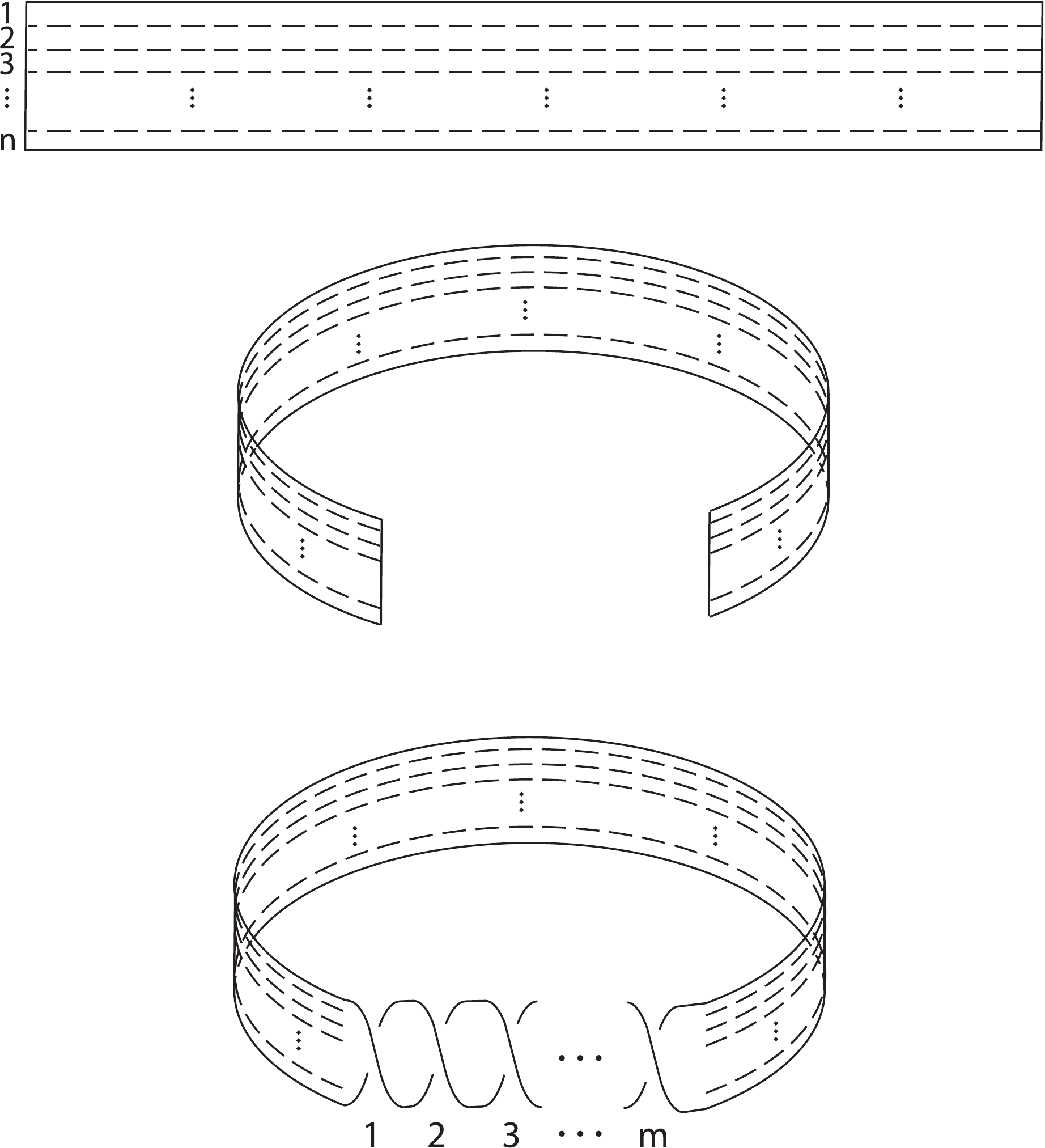}
\caption{\label{figpconst}%
The $P(m,n)$ paradromic ring: join the ends with $m$ half-twists and cut
along the dashed lines.}
\end{center}
\end{figure}

Let's review how we construct a $P(m,n)$ 
paradromic ring (see Figure~\ref{figpconst}).
Draw lines on a  strip of paper that divide it into $n$ strips. Connect the two loose ends with $m$ half twists and then cut along the lines. Finally, we replace each resulting loop of paper, whose width is $1/n$ that
of the original strip, with the curve that
runs along its midline, $1/2n$ from its edges. 
We assume $m$ is a non-negative integer and $n$ is positive. 

\textbf{FIGURE 11 GOES NEAR HERE}.

\begin{figure}[p]
\begin{center}
\includegraphics[scale=.35]{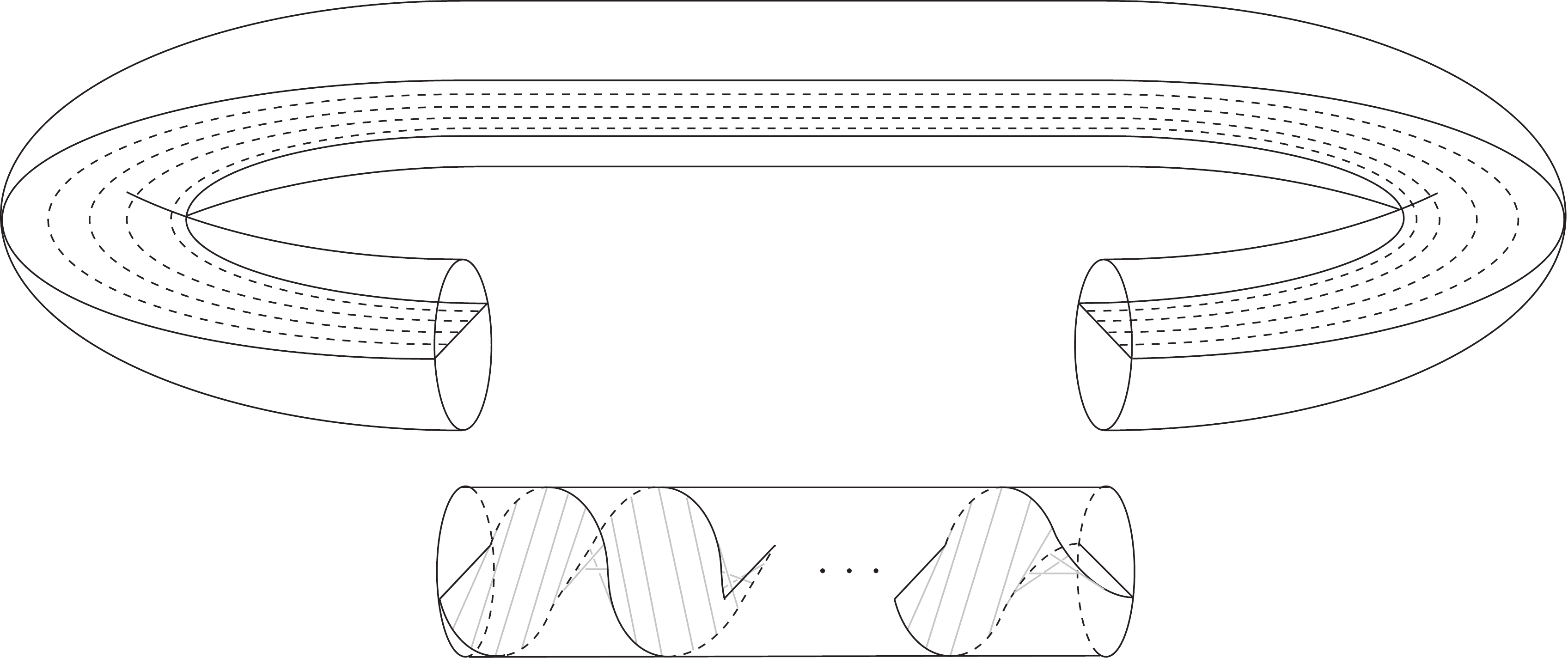}
\caption{\label{figpint}%
Isolate the twists in a cylinder, $C_t$. Outside the cylinder, the strip lies between the inner and outer equators on the torus.}
\end{center}
\end{figure}

To illustrate the connection with torus links, we place
our strip of paper inside a torus (see Figure~\ref{figpint}). 
We will group the $m$ half twists together (compare with the
$P(m,2)$ diagrams at the bottom of Figure~\ref{FigP352braid})
and then connect them up with a flat strip that joins the two ends of
the twisted region. In other words, we collect the half twists 
inside a cylinder that we'll call $C_t$ ($t$ for twist). Outside the cylinder, the strip of paper will lie between 
concentric circles that we call the {\em equators}.

\textbf{FIGURE 12 GOES NEAR HERE}.

\begin{figure}[p]
\begin{center}
\includegraphics[scale=.9]{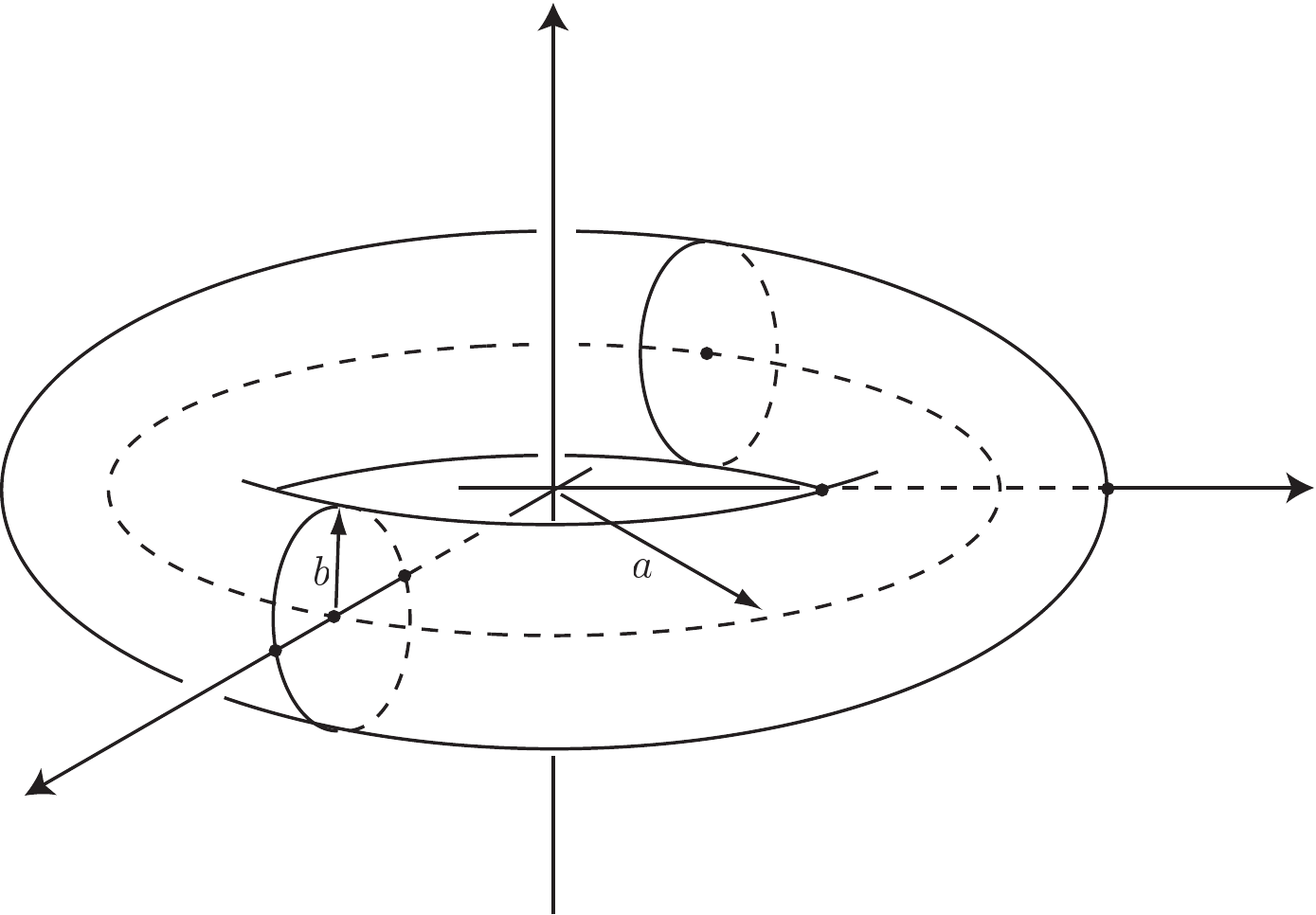}
\caption{\label{figtorpar}%
An embedding of the torus in $\R^3$. The $z$-axis is an axis of rotational symmetry. The $xy$-plane is fixed by a reflection.}
\end{center}
\end{figure}

For convenience in defining equators, the core, and other nomenclature, we 
situate the torus in $\R^3$ as in 
Figure~\ref{figtorpar}. The 
$z$-axis is an axis of rotational symmetry and
the $xy$-plane is fixed by a reflection. Let $a$ and $b$ be 
the radii shown in the figure. The core, then, is the circle
in the $xy$-plane of radius $a$ centered at the origin.
The $xy$-plane  intersects the torus in two concentric circles (of radius $a-b$ and $a+b$) that we call 
the {\em inner} and {\em outer equators.} A {\em longitude} is any closed curve on the torus that is parallel to the equators and loops once around the $z$-axis. For example, planes of the form $z = c$ where $|c| < b$ will intersect the torus in two longitudes. The plane $z = b$ intersects
the torus in a single longitude, the {\em top longitude}, that runs along the top of the torus. The equators are also examples of longitudes.
A {\em meridian} is any simple closed curve that
intersects each longitude once and also bounds a disk inside the torus.
Planes of the form $y = kx$, for example, intersect the torus in two meridia, each being a circle of radius b.

The {\em $T(u,v)$ torus link} is a link of $\mbox{GCD}(u,v)$ components 
that we can arrange on the torus so that it
intersects each longitude $u$ times and each meridian $v$ times.
As mentioned in
Section~1, when we speak of a link, an embedding of circles in 
three space, we are allowed to move the circles around in space freely
so long as the curves do not pass through one another.
Such a link is a torus link if, among these different embeddings, 
there is one that lies flat on a torus without the curve crossing through
itself. 
For example, in Figure~\ref{figtorex}, the trefoil is $T(3,2)$, the pentafoil is  $T(5,2)$, and the Hopf link is  $T(2,2)$. We have starred the intersections with the outer equator, which is a longitude.

We are now ready to prove Theorem~1: either a paradromic ring is a torus link, 
or else it is a torus link together with an additional component along the core 
of the torus. We denote the second case by $T(u,v) \cup C$. 
Figure~\ref{figcore} shows, for example, that the $P(3,3)$ paradromic ring
is $T(3,2) \cup C$. 

\begin{thm}
Let $m \geq 0$ and $n > 0$ be integers.
If $n=1$, 
$P(m,1) = T(0,1)$; if $n>1$, then
$$
P(m,n) = \left\{ \begin{array}{ll}
T(\frac12mn,n)  & \mbox{ if } mn \mbox{ is even,} \\
\\
T(\frac12m(n-1),n-1) \cup C & \mbox{ if } mn \mbox{ is odd. }
\end{array} \right.
$$
\end{thm}
\noindent%

\bigskip

Below we sketch an argument that is largely a proof by pictures. This is a
perfectly respectable technique used by professional topologists the 
world over. We could, if needed, replace it with an `analytic' proof
that doesn't rely on pictures, but that would be very tedious and less insightful. 

Still, if the idea of a proof by pictures is not to your
taste, we encourage you to accept the theorem for the 
sake of argument and skip ahead to Section~3 where
linear algebra again comes to the fore.

\bigskip

\textbf{FIGURE 13 GOES NEAR HERE}.

\begin{figure}[p]
\begin{center}
\includegraphics[scale=.35]{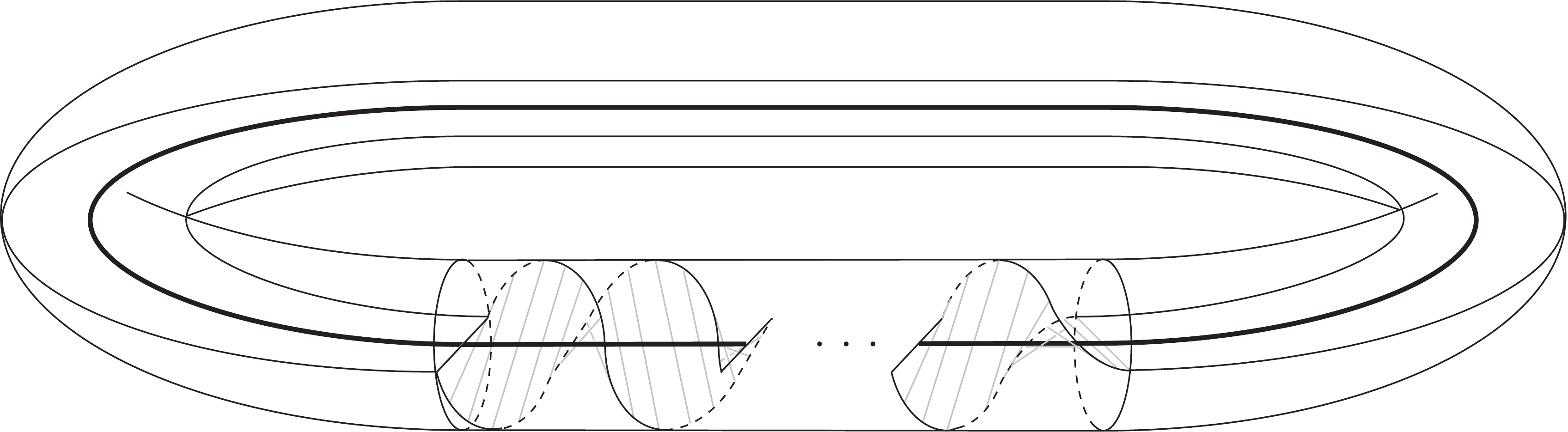}
\caption{\label{fign1}%
If $n=1$ the midline (bold) follows the core of the torus.}
\end{center}
\end{figure}

\Pf (sketch)
If $n=1$, we do not cut the strip of paper at all; it consists of a single
loop whose midline follows the core of the torus, see Figure~\ref{fign1}. Moving the core straight up in the $z$-direction to follow the top longitude, 
we see that $P(m,1) = T(0,1)$. In other words, as a knot, the core is equivalent to any longitude since we can move it in space to follow that longitude.

\textbf{FIGURE 14 GOES NEAR HERE}.

\begin{figure}[p]
\begin{center}
\includegraphics[scale=.35]{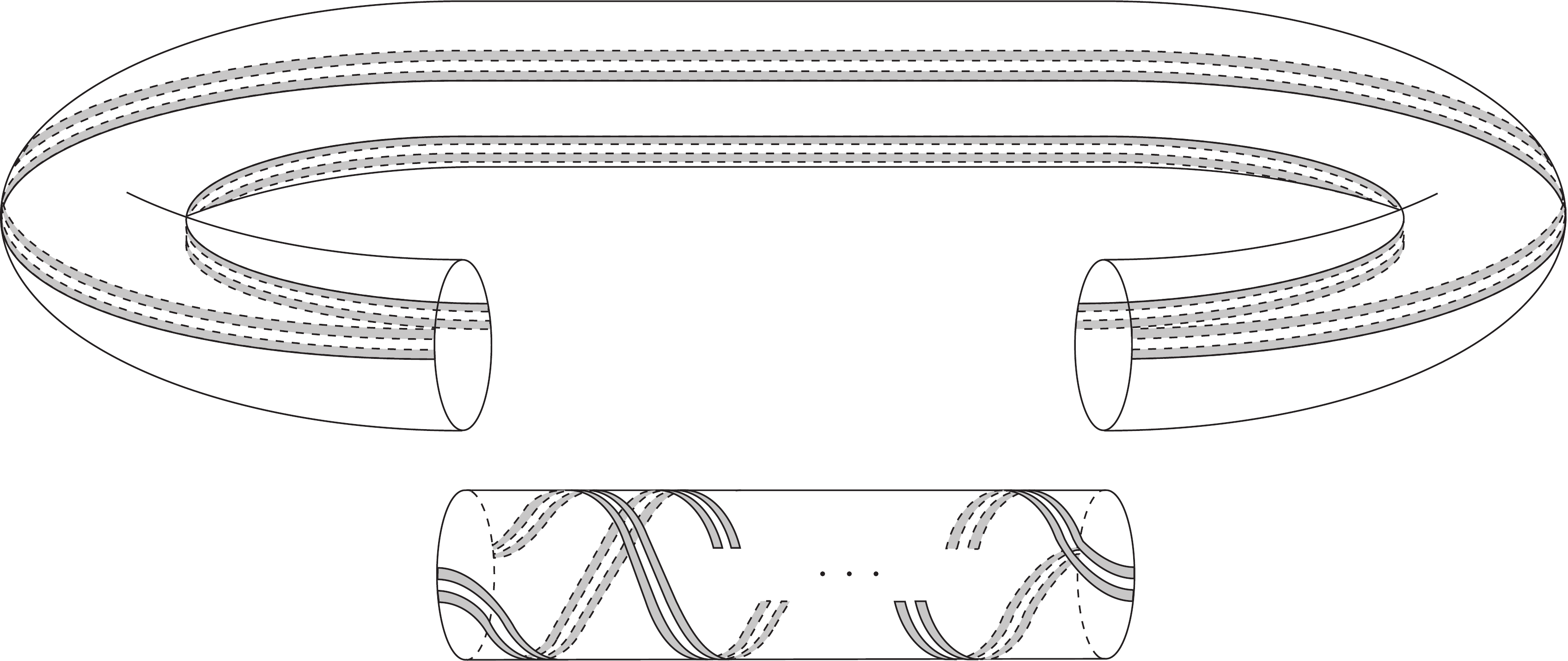}
\caption{\label{figpusheven}%
If $n$ is even, after halving, the $n$-sections can be pushed into the torus. Here, $n = 4$.}
\end{center}
\end{figure}

When $n > 1$, we 
place our twisted strip of paper inside a torus, as in Figure~\ref{figpint},
with all twists gathered in the cylinder $C_{t}$ ($t$ for twist).
If $n$ is even, then one of the dashed lines of Figure~\ref{figpconst}
will run right down the center of the strip. Cutting along this line bisects the strip
and allows us to lay the bisected strip flat on the torus. 
(We are taking advantage of the idea that we are free to move a link
around in space so long as we do not pass it through itself.)
Outside of $C_{t}$,
we can think of the strip's two halves as two narrow bands, one near the inner equator and one near the outer equator (see Figure~\ref{figpusheven}).

After cutting the strip into its $n$ sections, we will have a 
collection of thin strips on the torus, half grouped around the 
inner equator and half around the outer equator.
Outside of $C_{t}$, this collection of strips
cross a meridian $n$ times, with $n/2$ intersections near each of the 
two equators. On the other hand, the strips will 
cross a longitude $mn/2$ times. For example, the top longitude intersects the rings only in $C_{t}$, and there we have $n/2$ crossings for each half twist. Thus, we have a $T(mn/2,n)$ torus link.

\textbf{FIGURE 15 GOES NEAR HERE}.

\begin{figure}[p]
\begin{center}
\includegraphics[scale=.35]{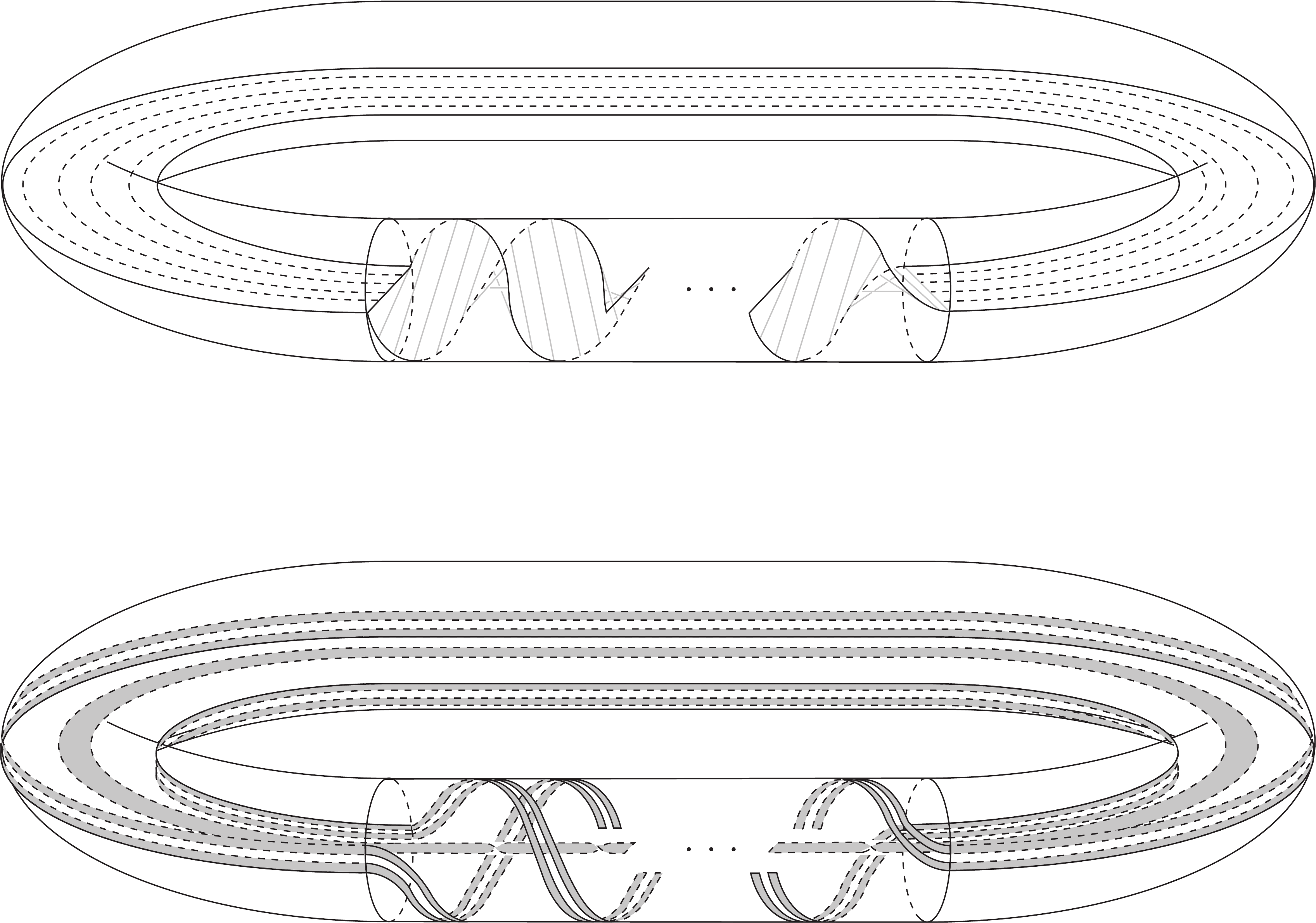}
\caption{\label{figpushodd}%
If $n$ is odd, going from top to bottom, we leave the central strip at the core and push the remaining $n-1$ sections onto the torus. Here, $n = 5$.}
\end{center}
\end{figure}

If $n$ is odd, by leaving the central strip at the core of the torus,
we can again place the remaining $n-1$ sections onto the torus with $(n-1)/2$ strips near each of the two equators, see Figure~\ref{figpushodd}.
In addition to the core, we are left with strips on the torus
that cross each meridian $n-1$ times while meeting a longitude
$m(n-1)/2$ times, resulting in $T(m(n-1)/2,n-1) \cup C$.

Finally, if $n$ is odd and $m$ is even, we can also move the strip at the core onto the torus, making a torus link. For example, move the core to follow the top longitude
outside of $C_{t}$. If we continue the curve into $C_{t}$ 
starting at the top of the cylinder at left, then after $m$ (an even number) of half twists, it will have returned to the top when we reach 
the right end of $C_{t}$ so that we can close the curve.
Compared to $T(m(n-1)/2, n-1)$, this adds
an extra intersection with each meridian and $m/2$ intersections with each
longitude. This is the $T(mn/2,n)$ torus link. \qed

\section{Paradromic rings resist coloring}

We are now ready to classify the colorability of the paradromic rings.
We break the argument into two cases, as in Theorem~1: 
paradromic rings that are torus links, and those that are not. 

\textbf{FIGURE 16 GOES NEAR HERE}.

\begin{figure}[p]
\begin{center}
\includegraphics[scale=.6]{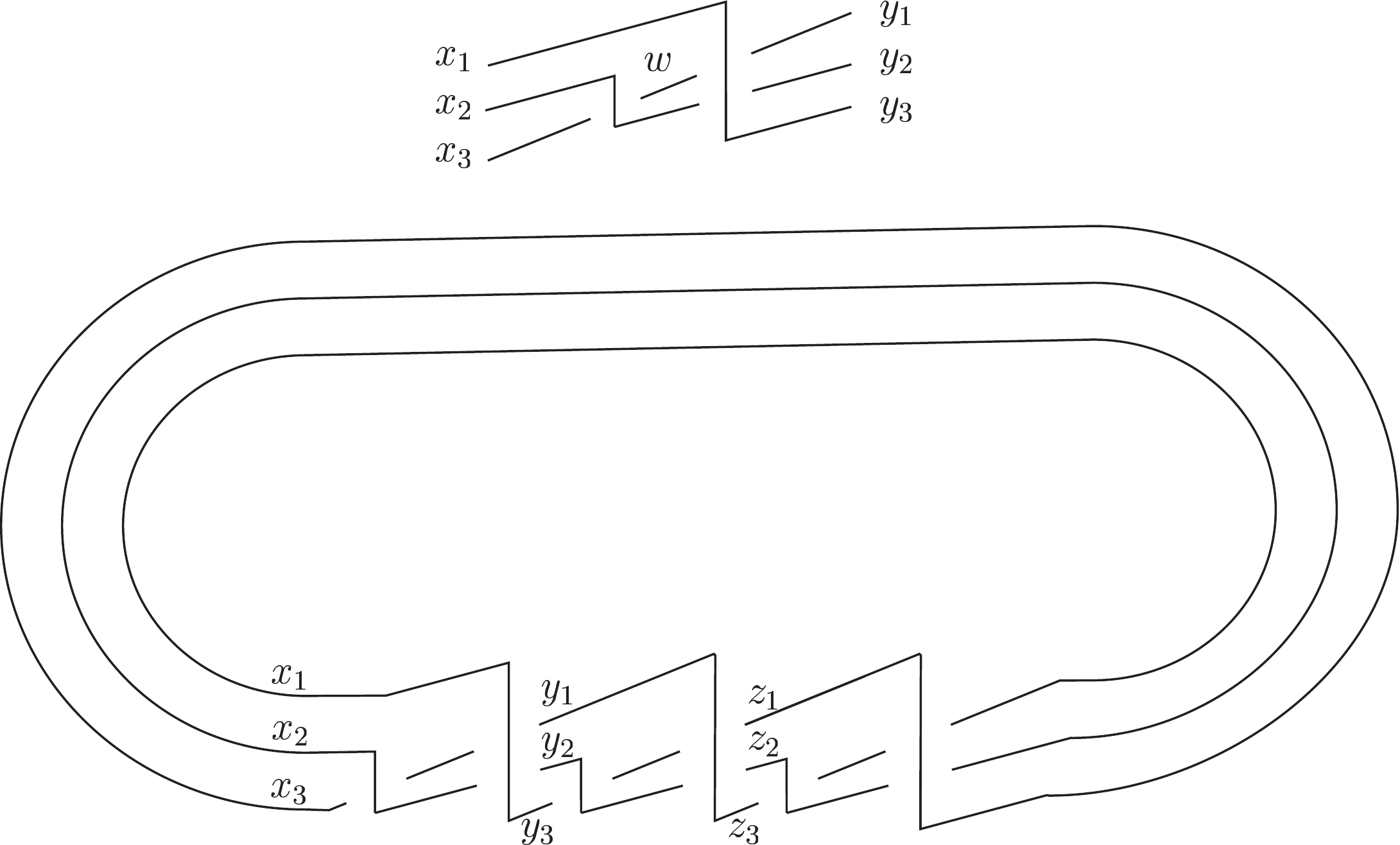}
\caption{
\label{figP33}%
$P(3,3)$ is formed by repeating the pattern three times.}
\end{center}
\end{figure}

We begin with those that are not, in other words, the $P(m,n)$ where $mn$ is odd and $n > 1$. 
The $P(m,2)$ torus links of Section~1 illustrate
our approach. As a further example, let's color $P(3,3)$, which is 
not a torus link (see Figure~\ref{figcore}).
Figure~\ref{figP33} shows how to construct this link
by repeating the pattern at top three times. 
Color the arcs entering the pattern at left with $x = (x_1, x_2, x_3)$. Then a matrix equation
determines the colors 
$y = (y_1, y_2, y_3)$ leaving at right:
$y \equiv S_3 x$. 

Let's find the matrix $S_3$. Referring to the pattern at the top of Figure~\ref{figP33}, there are two crossings involving $x_1$, both with $x_1$ as the overarc.
In the lower one, condition 2 for $p$-colorability yields
$2x_1  \equiv x_2 + y_2 \pmod{p} \Rightarrow
y_2 \equiv 2x_1 - x_2 \pmod{p}$. At the upper crossing, we have 
$2x_1 \equiv w + y_1  \pmod{p} \Rightarrow 
y_1 \equiv 2x_1  - w \pmod{p}$. The third crossing in the pattern shows how to write $w$ in
terms of $x_2$ and $x_3$: $2x_2  \equiv x_3 + w \pmod{p} \Rightarrow w \equiv 2x_2 - x_3 \pmod{p} $. Thus, we have the following system of equations modulo $p$:
\begin{eqnarray*}
2x_1 - (2x_2 - x_3) & \equiv & y_1 \\
2x_1 - x_2 & \equiv & y_2 \\
x_1 & \equiv & y_3
\end{eqnarray*}
with coefficient matrix
$$
S_3 = \left(
\begin{array}{rrr}
2 & -2 & 1 \\
2 & -1 & 0 \\
1 & 0 & 0 
\end{array}
\right) .
$$

Similarly, $(z_1,z_2,z_3) = z \equiv S_3 y \pmod{p}$. Following the arcs around
the top of the link, we see that $x \equiv S_3 z \pmod{p}$. This means a $p$-coloring of $P(3,3)$
corresponds to a vector $x$ such that $x \equiv S_3^3 x \pmod{p}$. 
In other words, we want an eigenvector of $S_3^3$ modulo $p$ with eigenvalue one. 

The characteristic polynomial of $S_3^3$ is $\dt{S_3^3 - \lambda I} = 
-(\lambda -1)(\lambda^2+1)$. As long as $p \neq 2$, the $\lambda = 1$ eigenspace has dimension one and the
only eigenvectors are the constant vectors, $(c,c,c)$.
Recall that a constant vector means all arcs
in the diagram have color $c$, in violation of condition 1 for
$p$-coloring. Therefore, when $p \neq 2$, $P(3,3)$ is not $p$-colorable. 
On the other hand, as $P(3,3)$ has
two components, it is $2$-colorable. For example, we could color the core
$0$ and the trefoil component $1$. Thus,  
$P(3,3)$ is nearly invisible. It is $p$-colorable only for the prime
$p = 2$. 

As the following theorem shows, this is true of all the
paradromic rings that are not torus links. We began our study expecting 
that $p$-colorability would be an interesting way to distinguish among
these rings. Instead it turns out that they are all nearly invisible.

\begin{thm} 
\label{thmnear}
If $m$ and $n$ are positive odd integers with $n > 1$, then the paradromic ring $P(m,n)$ is nearly invisible.
\end{thm}

\textbf{FIGURES 17 AND 18 GO NEAR HERE}.

\begin{figure}[p]
\begin{center}
\includegraphics[scale=.7]{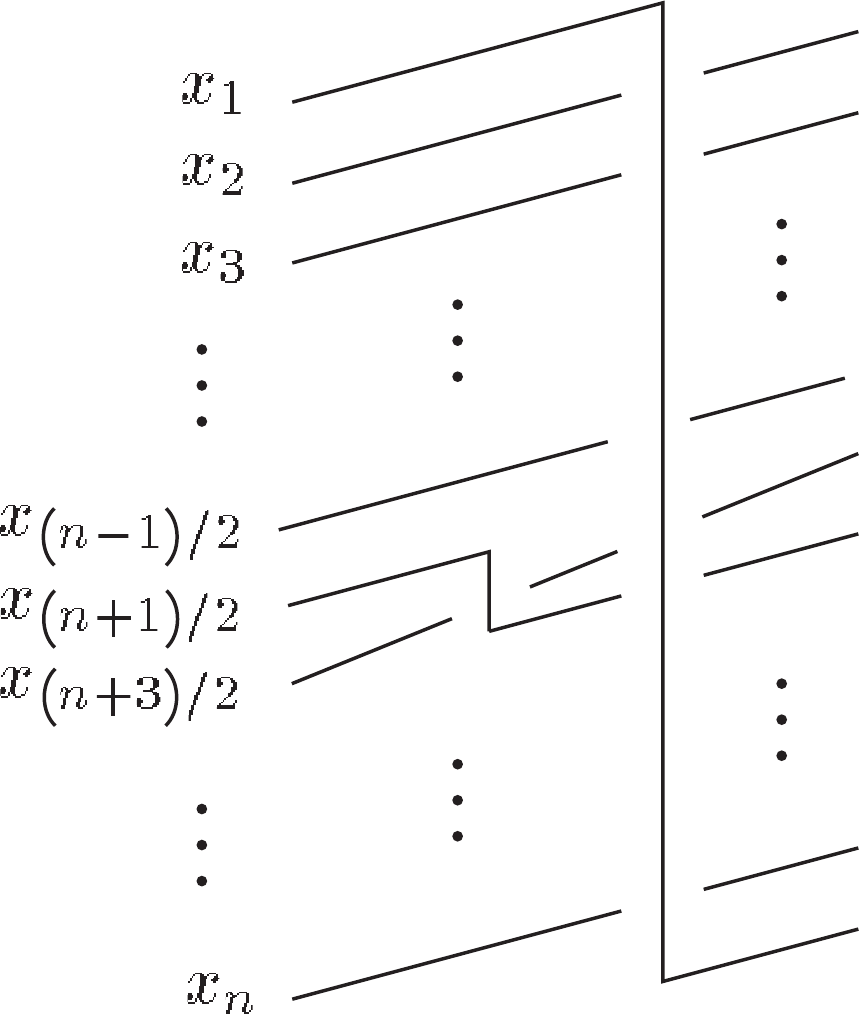}
\caption{\label{figApattn}%
A pattern on $n$ arcs (where $n > 1$ is odd). }
\end{center}
\end{figure}

\begin{figure}[p]
\begin{center}
\includegraphics[scale=.5]{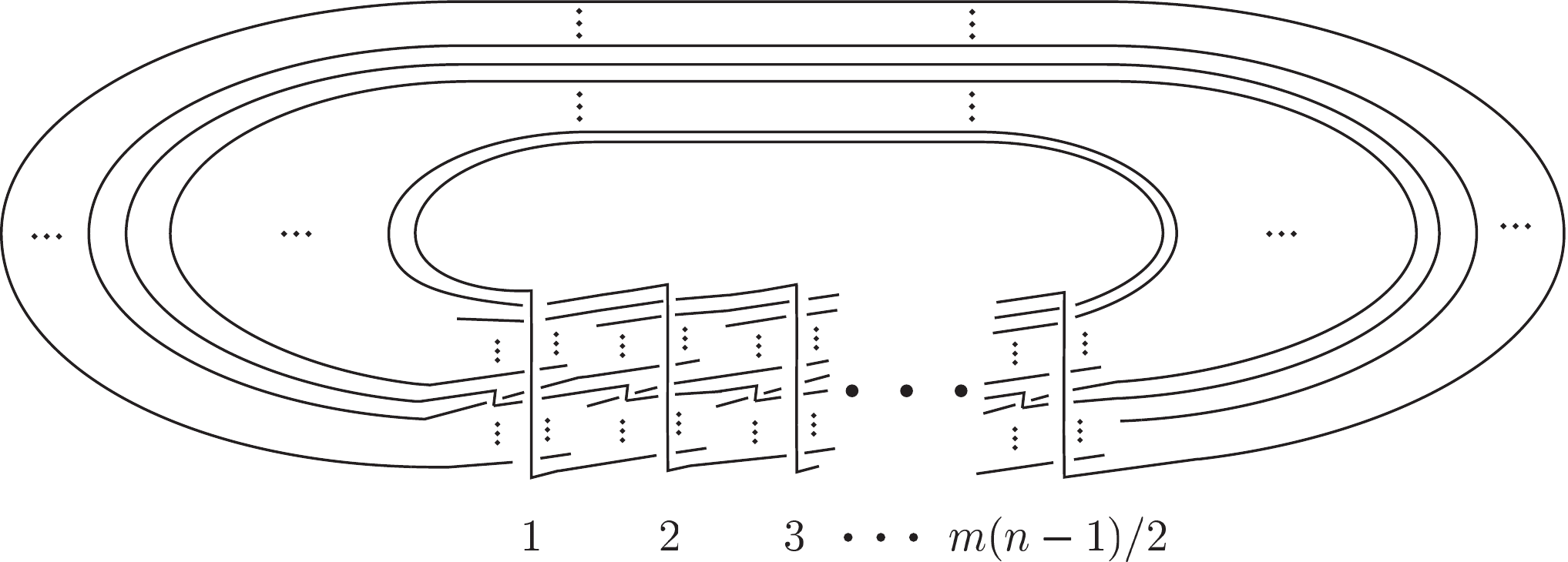}
\caption{\label{figApattn2}%
Repeating the pattern on $n$ arcs $m(n-1)/2$ times forms $P(m,n)$.}
\end{center}
\end{figure}

Before proving the theorem, we will describe the matrix $S_n$ that generalizes $S_3$ for $n$ odd. 
Let $m$ and $n$ be positive odd integers. 
We represent $P(m,n)$ as in Figure~\ref{figApattn2}, as suggested by 
our analysis of $P(m,2)$ and $P(3,3)$.
That is, $P(m,n)$ consists of $m (n-1)/2$ repetitions
of the pattern in Figure~\ref{figApattn} joined up in a ring. 
This figure gives us the matrix
$$
S_n = \left(
\begin{array}{rrrrrrrrr} 
2 & -1 &  0 & \ldots &  0 & 0 & \ldots &  0 &  0 \\
2 &  0 & -1 & \ldots &  0 & 0 & \ldots &  0 &  0 \\
\vdots & \vdots & \ddots & \ddots & \vdots & \vdots & \ldots & \vdots & \vdots \\
2 &  0 &  0 & \ddots & -2 & 1 & \ldots &  0 &  0 \\
2 &  0 &  0 & \ldots & -1 & 0 & \ldots &  0 &  0 \\
\vdots & \vdots & \vdots & \ldots & \vdots & \ddots & \ddots & \vdots & \vdots \\
2 &  0 &  0 & \ldots &  0 & 0 & \ddots & -1 &  0 \\
2 &  0 &  0 & \ldots &  0 & 0 & \ldots &  0 & -1 \\
1 &  0 &  0 & \ldots &  0 & 0 & \ldots &  0 &  0 
\end{array} \right) .
$$
If $x = (x_1, \ldots, x_n)$ are the colors of the arcs entering the pattern of Figure~\ref{figApattn} at the left, 
then the outgoing arcs at right are $S_n x$ modulo $p$. 
Note that, outside of a $2 \times 2$ block, $S_n$ has $-1$'s on the superdiagonal and a first column that is all $2$'s but for a $1$ in the last row. The $2 \times 2$ matrix, 
$$
\left( \begin{array}{rr}
-2 & 1 \\
-1 & 0 
\end{array} \right),
$$
that breaks up the pattern is 
in rows $(n-1)/2$ and $(n+1)/2$ and columns $(n+1)/2$ and $(n+3)/2$
(recall that $n > 1$ is odd) and is due to the short $w$ arc
in the middle of the pattern.
The $S_n$ matrix has a surprisingly simple characteristic polynomial.

\begin{lemma} 
\label{lem41}
Let $n > 1$ be an odd integer.
The characteristic polynomial of $S_n$ is $f_n(\lambda) = -(\lambda-1)(\lambda^{n-1} +1)$.
\end{lemma}

\Pf
Since $S_n$ follows a regular pattern except for columns $(n+1)/2$ and $(n+3)/2$, we will make 
expansions along those columns to recover more symmetric matrices. 
Expanding along column
$(n+3)/2$,  $f_n(\lambda) = \dt{S_n - \lambda I}  = \dt{A_{n-1}} - \lambda \dt{B_{n-1}}$ where $A_{n-1}$ and
$B_{n-1}$ are $(n-1) \times (n-1)$ minors. 
Column $(n+1)/2$ then shows $$\dt{B_{n-1}} = 2 \dt{C_{n-2}} - (\lambda +1) \dt{D_{n-2}}.$$ 
Below, we argue
\begin{eqnarray*}
\dt{A_{n-1}} & = & (\lambda+1)  - 2 \lambda (1 - (- \lambda)^{\frac{n-3}{2}}) \\
\dt{C_{n-2}} & = & 2 (-\lambda)^{\frac{n-3}{2}}  \mbox{, and }\\
\dt{D_{n-2}} & = & -\lambda^{\frac{n-3}{2}} \left( \frac{\lambda^{\frac{n+1}{2}}-  \lambda^{\frac{n-1}{2}}+ 2 (-1)^{\frac{n-1}{2}}}{\lambda+1} 
\right).
\end{eqnarray*}

Then, we have 
\begin{eqnarray*}
f_n (\lambda) & = & \dt{A_{n-1}} - \lambda \dt{B_{n-1}} \\
& = & \dt{A_{n-1}} - \lambda (2 \dt{C_{n-2}} - (\lambda +1) \dt{D_{n-2}}) \\
& = & -(\lambda-1)(\lambda^{n-1} +1).
\end{eqnarray*}

Let's verify the formulas for the determinants of $A_{n-1}$, $C_{n-2}$, and $D_{n-2}$.
After appropriate column and row expansions (Start with column $(n+1)/2$.) we deduce
$\dt{A_{n-1}} = (\lambda+1)(1-\lambda \dt{\bar{A}_{(n-1)/2}})$
where $\bar{A}_k$ is the $k \times k$ matrix
$$\bar{A}_k =
\left(
\begin{array}{crrrrrrr} 
2- \lambda  & -1 &  0 &  0 & \ldots &  0 &  0 &  0 \\
2 &  0 & -1 &  0 & \ldots &  0 &  0 &  0 \\
2 & 0 & -\lambda & -1 & \ldots & 0 & 0 & 0 \\
\vdots & \vdots & \vdots & \ddots & \ddots & \vdots & \vdots & \vdots \\
2 &  0 &  0 & 0 & \ddots &  -1 & 0 &  0 \\
2 &  0 &  0 & 0 & \ldots &  -\lambda & -1 &  0 \\
2 &  0 &  0 & 0 & \ldots &  0 &  -\lambda & -1 \\
2 &  0 &  0 & 0 & \ldots &  0 &  0 &  -\lambda 
\end{array} \right). 
$$
Expanding along the last row, we find $\dt{\bar{A}_k} = 2 - \lambda \dt{\bar{A}_{k-1}}$. Solving the recurrence relation, we have
$$\dt{\bar{A}_k} = 2 \left(  \frac{1 - (- \lambda)^{k-1}}{1+\lambda} \right),$$ as required.

For $C_{n-2}$, the $(n-1)/2$ row is zero but for a 2 at the beginning of the row. Expanding along that row, 
we uncover a minor that is a block diagonal matrix. The top left block is lower triangular with determinant $(-1)^\frac{n-3}{2}$ and the
bottom right block is upper triangular with determinant $(-\lambda)^{\frac{n-3}{2}}$. The sign of the
determinant depends on the parity of $(n-1)/2$, the row along which we expand.

Much like $A_{n-1}$, we express $\dt{D_{n-2}}$ in terms of a smaller, more symmetric matrix:
$\dt{D_{n-2}} = (-\lambda)^{\frac{n-1}{2}} \dt{\bar{D}_{\frac{n+1}{2}}}$ where
$$\bar{D}_k =
\left(
\begin{array}{crrrrrrr} 
2- \lambda  & -1 &  0 &  0 & \ldots &  0 &  0 &  0 \\
2 &  -\lambda & -1 &  0 & \ldots &  0 &  0 &  0 \\
2 & 0 & -\lambda & -1 & \ldots & 0 & 0 & 0 \\
\vdots & \vdots & \vdots & \ddots & \ddots & \vdots & \vdots & \vdots \\
2 &  0 &  0 & 0 & \ddots &  -1 & 0 &  0 \\
2 &  0 &  0 & 0 & \ldots &  -\lambda & -1 &  0 \\
2 &  0 &  0 & 0 & \ldots &  0 &  -\lambda & -1 \\
2 &  0 &  0 & 0 & \ldots &  0 &  0 &  -\lambda 
\end{array} \right). 
$$
Again, $\dt{\bar{D}_k} = 2 - \lambda \dt{\bar{D}_{k-1}}$, and solving
the recurrence, yields the formula for $\dt{D_{n-1}}$.
\qed

\medskip

\Pf (of Theorem~2)
Let $u = m(n-1)/2$ and let $p$ be an odd prime. 
Colorings of $P(m,n)$ are $\lambda = 1$ eigenvectors of $S_n^u$ modulo $p$.
We will show that $\lambda = 1$ is a simple root of the characteristic polynomial of $S_n^u$. 
This means the only eigenvectors are the constant vectors $(c,c,c, \ldots, c)$
and there are no valid colorings when $p$ is odd.
Since $P(m,n)$ has at least two components, the core and a torus link, it is $2$-colorable.
This shows that $2$ is the only prime coloring and $P(m,n)$ is nearly invisible.

Let's see why $\lambda = 1$ is a simple root when $p$ is odd.
Let $F$ be the characteristic polynomial of $S_n^u$. The roots of
$F$ are the $u$th powers of the roots of $f_n$, the characteristic polynomial of $S_n$. By Lemma~1, $1$ is a root of $f_n$ and hence of $F$.

We must argue that no other root of $F$ is equal to $1$.
That is, if $\zeta$ is a root of the second factor of $f_n$, $x^{n-1} + 1$,
we must show $\zeta^u \not \equiv 1 \pmod{p}$.

Let $\zeta$ be a root of $x^{n-1} + 1$. Then $\zeta^{n-1} \equiv -1$. Suppose, for a contradiction, that $\zeta^u \equiv 1 \pmod{p}$. Now, 
since $m(n-1)  = 2 u $,
\begin{eqnarray*}
\zeta^{m(n-1)} \equiv \zeta^{2 u}
& \Rightarrow & (-1)^m \equiv 1^2 \\
& \Rightarrow & -1 \equiv 1 \pmod{p}, 
\end{eqnarray*}
which is absurd since $p$ is not $2$. 

The contradiction shows
that the roots of $x^{n-1}+1$ do not lead to additional
occurences of $1$ as a root of the characteristic polynomial $F$ of
$S_n^u$. Therefore, $S_n^u$ has no non-constant eigenvectors with
eigenvalue one and $P(m,n)$ is not $p$-colorable for any odd prime
$p$.
\qed

The paradromic rings that are torus links include infinite families of
rainbow rings and nearly invisible links:
\begin{thm} 
\label{thmrainbow}
Let $n>1$ and $m \geq 0$ be integers such that $mn$ is even.
Then the torus link  $T = T(\frac12mn,n)$ is a rainbow ring unless one of the following occurs:
\begin{itemize}
\item $n$ and $\frac12m$ are both odd, in which case $T$ is nearly invisible
\item $n = 2$, in which case $T$ is $p$-colorable if and only if $p$ divides $m$
\item $n=4$ and $m$ is odd, in which case $T$ is $p$-colorable if and only if $p$ divides $2m$.
\end{itemize}
\end{thm}

On the other hand, many of the knots in the family are invisible:
when $n = 1$, $P(m,1)$ is just a circle whose determinant is one.

We omit the proof 
of Theorem~\ref{thmrainbow}
for a couple of reasons. First, we expect that an inspired reader 
is capable of completing the proof,
just as the REU team did during the summer.
In particular, Section~1 above includes
the argument for $P(m,2)$ (that is, the case where $n = 2$). 

Second, we want to take the chance to recommend additional 
reading that leads to a more direct approach in the case of torus links.
The colorability of torus {\bf knots} has already been determined by other researchers including Bryan~\cite{B}, and Breiland, Oesper, and Taalman~\cite{BOT}:

\begin{thm}[\cite{B, BOT}] 
\label{thmknots}
Let $u,v$ be positive integers with $\mbox{GCD}(u,v) =1$. The torus knot $T(u,v)$ is $p$-colorable if and only if either $u$ is even and $p$ divides $v$ or else $v$ is even and $p$ divides $u$.
\end{thm}

Indeed, it was Bryan's analysis that inspired us to attempt a similar
argument for paradromic rings.

We have already recommended
Adams's {\em The Knot Book}~\cite{A} and Livingston's {\em Knot Theory}~\cite{L} 
as nice introductions to $p$-coloring, including the proof that it is a link invariant.
Murasugi's {\em Knot Theory \& Its Applications}~\cite{M} is at a slightly more advanced level and includes a thorough introduction to 
the idea of the determinant of a link, $\dt{L}$, and how to calculate it.
As you will read there, $\dt{L}$ is indeed the determinant of a matrix,
although not the matrices $S_n$ and $T$ discussed in this paper. 
Making use of that matrix, Murasugi shows that 
the determinant of a  torus link $L  = T(u,v)$
is given by $\dt{T(u,v)} = | \Delta(-1) |$ where, up to a multiple of $x$,
$$\Delta(x) = \frac{(1-x)(1-x^{\frac{uv}{d}})^d}{(1-x^u)(1-x^v)},$$
with $d = \mbox{GCD}(u,v)$.  
Recalling that a link $L$ is $p$-colorable if and only 
$p$ divides $\dt{L}$, the formula gives a direct way to prove Theorem~3.
In particular, when $n \geq 5$, the GCD $d$ is at least 3, which means that terms of the form $1-x^{2k}$ survive in the
numerator so that $|\Delta(-1)| = 0$ (provided $n$ and $m/2$ are not both odd).

\paragraph*{Acknowledgements} This paper grew out of a 2005 REUT at
CSU, Chico that was supported in part by
NSF REU Award 0354174 and by the MAA's NREUP program 
with funding from the NSF, NSA, and Moody's. The first three authors were undergraduates at the time while Dan Sours is a high school teacher.
We are grateful to Yuichi Handa, Ramin Naimi, Neil Portnoy, Robin Soloway, and John Thoo for helpful comments on early versions of this paper.
Additional funding came from CSU, Chico's CELT as part of a 2015 Faculty Learning Community. We thank 
Chris Fosen, Greg Cootsona, and the other FLC participants for fruitful discussions about the exposition.

\end{document}